\documentstyle[12pt]{article}
\newtheorem{thm}{Theorem}[section]

\newtheorem{pro}[thm]{Proposition}

\newtheorem{rem}[thm]{Remark}

\def\b1{\mbox{\boldmath $1$}}

\begin{document}
\begin{center}
{\large\bf Some Diffusion Processes Associated With Two Parameter
Poisson-Dirichlet Distribution and Dirichlet Process}
\end{center}
\vskip 0.5cm
\begin{center} {Shui Feng and Wei
Sun} \vskip 0.1cm {\it McMaster University and Concordia
University} \vskip 0.5cm
\begin{minipage}{120mm}
{\parindent 6mm {\small \noindent The two parameter
Poisson-Dirichlet distribution $PD(\alpha,\theta)$ is the
distribution of an infinite dimensional random discrete probability.
It is a generalization of Kingman's Poisson-Dirichlet distribution.
The two parameter Dirichlet process $\Pi_{\alpha,\theta,\nu_0}$ is
the law of a pure atomic random measure with masses following the
two parameter Poisson-Dirichlet distribution. In this article we
focus on the construction and the properties of the infinite
dimensional symmetric diffusion processes with respective symmetric
measures $PD(\alpha,\theta)$ and $\Pi_{\alpha,\theta,\nu_0}$. The
methods used come from the theory of Dirichlet forms.}}
\end{minipage}
\end{center}
\vskip 0.3cm

\section[short title] {Introduction}
\setcounter{equation}{0}

 The Poisson-Dirichlet distribution $PD(\theta)$ was introduced by Kingman in \cite{Kingman75}
 to describe the distribution of gene frequencies in a large neutral
      population at a particular locus. The component $P_k(\theta)$ represents the
      proportion of the $k$-th most frequent allele. The Dirichlet process  $\Pi_{\theta,\nu_0}$ first appeared
      in \cite{Fer73} in the context of Bayesian statistics. It is a pure atomic random measure with masses distributed
      according to $PD(\theta)$.

In the context of population genetics, both the Poisson-Dirichlet distribution and the Dirichlet process appear
   as approximations to the equilibrium behavior of certain large populations evolving under the influence of mutation
   and random genetic drift. To be precise, let $C_b(S)$ be the set of bounded, continuous
functions on a locally compact, separable metric space $S$, ${\cal
M}_1(S)$ denote the space of all probability measures on $S$
equipped with the usual weak topology, and $\nu_0\in{\cal
M}_1(S)$. We consider the operator $B$ of the form
$$Bf(x) = \frac{\theta}{2}\int(f(y)-f(x))\nu_0(dy), f \in C_b(S).$$
Define $${\cal D} = \{u: u(\mu) = f(\langle \phi,\mu  \rangle),\,
f \in C^{\infty}_b({\mathbf{R}}), \phi \in C_b(S), \mu \in {\cal
M}_1(S) \},$$ where $C^{\infty}_b(\mathbf{R})$ denotes the set of
all bounded, infinitely differentiable functions on $\mathbf{R}$.
Then the Fleming-Viot process with neutral parent independent
mutation or the labeled infinitely-many-neutral-alleles model is a
pure atomic measure-valued Markov process with generator
 \[
 L u(\mu)  = \langle B\nabla u(\mu)(\cdot),\mu \rangle + \frac{f^{\prime\prime}(\langle \phi,\mu\rangle)}{2}
\langle \phi, \phi\rangle_{\mu}, u\in{\cal D},\] where
\begin{eqnarray*}
&& \nabla u(\mu)(x)=\delta u(\mu)/\delta \mu(x) = \lim_{\varepsilon \rightarrow 0+ }\varepsilon^{-1}
\{u((1-\varepsilon)\mu + \varepsilon \delta_x)-u(\mu)\},\\
&&\langle \phi,\psi \rangle_{\mu}=\langle \phi\psi, \mu\rangle-\langle \phi,\mu\rangle \langle \psi, \mu\rangle,
\end{eqnarray*} and
$\delta_x$ stands for the Dirac measure at $x \in S$. For compact
space $S$  and diffusive probability $\nu_0$, i.e., $\nu_0(x)=0$
for
      every $x$ in $S$, it is known
(cf. \cite{Et90}) that the labeled infinitely-many-neutral-alleles
model is reversible with reversible measure $\Pi_{\theta,\nu_0}$.

Introduce a map $\Phi$ from ${\cal M}_1(S)$ to the infinite
dimensional ordered simplex
\[
\nabla_{\infty} =\{(x_1,x_2,\ldots): x_1\geq x_2\geq \cdots \geq 0, \sum_{i=1}^{\infty}x_i=1 \}
\]
so that $\Phi(\mu)$ is the ordered masses of $\mu$. Then the labeled infinitely-many-neutral-alleles model is mapped
through $\Phi$ to another symmetric diffusion process, called the unlabeled infinitely-many-neutral-alleles model,
 on $\nabla_{\infty}$ with generator
\begin{equation}\label{A0}
A^0
=\frac{1}{2}\left\{\sum_{i,j=1}^{\infty}x_i(\delta_{ij}-x_j)\frac{\partial^2}{\partial
x_i\partial x_j} -\sum_{i=1}^{\infty}\theta x_i \frac{\partial
}{\partial x_i}\right\},
\end{equation}
defined on an appropriate domain. The symmetric measure of this
process is $PD(\theta)$.

  For any $0\leq \alpha <1$ and $\theta>-\alpha$, let $U_k,
      k=1,2,...$, be a sequence of independent
      random variables such that $U_k$ has $Beta(1-\alpha,\theta+ k\alpha)$ distribution.
      Set
      $$
      V^{\alpha,\theta}_1 = U_1,\  V^{\alpha,\theta}_n = (1-U_1)\cdots (1-U_{n-1})U_n,\  n \geq 2,
      $$
and let ${\bf P}(\alpha,\theta)=(\rho_1, \rho_2,...)$ denote
$(V^{\alpha,\theta}_1, V^{\alpha,\theta}_2,...)$
      in descending order. The distribution of $(V_1^{\alpha,\theta}, V_2^{\alpha,\theta}, \ldots)$ is called
      the two parameter GEM distribution. The law of ${\bf P}(\alpha,\theta)$ is called the two parameter Poisson-Dirichlet
      distribution, denoted by $PD(\alpha,\theta)$. For a locally compact, separable metric space $S$, and a sequence of i.i.d. $S$-valued random variables $\xi_k,
      k=1,2,...$ with common
      diffusive distribution $\nu_0$ on $S$, let
      \begin{equation}\label{diri1}
      \Xi_{\alpha,\theta, \nu_0}=\sum_{k=1}^{\infty}\rho_k\delta_{\xi_k}.
      \end{equation}
      The distribution of $\Xi_{\alpha,\theta, \nu_0}$, denoted by $Dirichlet(\theta,\alpha,\nu_0)$ or
      $\Pi_{\alpha,\theta,\nu_0}$, is called the two-parameter Dirichlet process. Clearly $PD(\theta)$ and $\Pi_{\theta,\nu_0}$ correspond to $\alpha=0$ in  $PD(\alpha,\theta)$
      and $\Pi_{\alpha,\theta,\nu_0}$, respectively.

As was indicated in \cite{PY} and the references therein, the two parameter Poisson-Dirichlet distribution and
 Dirichlet process are natural generalizations of their one parameter counterparts and possess many similar structures
including the urn construction, GEM representation, sampling
formula, etc. The cases of $\theta=0, \alpha$ are associated with
distributions of the lengths of excursions of Bessel processes and
Bessel bridge, respectively. It is thus natural to investigate the
two parameter generalizations of the labeled and unlabeled
infinitely-many-neutral-alleles models. One would hope that these
dynamical models will enhance our understanding of the two
parameter distributions.

 Several papers have appeared recently discussing the stochastic dynamics associated with the two parameter distributions.
 A symmetric diffusion process appears in \cite{fenwang07}, where the symmetric measure is the GEM distribution.
 An infinite dimensional diffusion process is constructed in \cite{P} generalizing the unlabeled infinitely-many-neutral-alleles model
 to the two-parameter setting. In \cite{berj08}, $PD(\alpha,\theta)$ is shown to be the unique reversible measure
 of a continuous time Markov chain constructed through
 an exchangeable fragmentation coalescence process. But it is still an open problem
 to construct the two parameter measure-valued process generalizing the Fleming-Viot
 process with parent independent mutation.

 In this article, we will consider two diffusion processes that are analogous to the
unlabeled and labeled infinitely-many-neutral-alleles models.
 In Section 2, an unlabeled two parameter infinitely-many-neutral-alleles
diffusion model is constructed via the classical gradient Dirichlet
forms. This process is shown to coincide with the process
constructed in \cite{P}. Besides establishing the existence and
uniqueness of the process, we also obtain results on the sample path
properties, the large deviations for occupation time process,
 and the model with interactive selection.
The construction of the labeled infinitely-many-neutral-alleles
diffusion model turns out to be much harder. Here the evolution of
the system involves both the masses and the locations. Note that
the one parameter model with finite many types is the
Wright-Fisher diffusion, and the partition property of the
infinite type model makes it possible for the finite dimensional
approximation. However, in the two parameter setting, the finite
type model itself is already a challenging problem not to mention
the loss of the partition property.  In Section 3, we construct a
general bilinear form that, if closable, will generate the needed
diffusion process. If the type space contains only two elements or
the type space is general but $\alpha=-\kappa$ and
$\theta=m\kappa$ for some $\kappa>0$ and integer $m \geq 2$, then
the above bilinear form is closable and a symmetric diffusion can
be constructed accordingly. The closability problem in the general
case boils down to the establishment of boundedness of a linear
functional. An auxiliary result is enclosed at the end of the
article to demonstrate the difficulty involved in establishing the
boundedness. If the bilinear form is indeed not closable, then its
relaxation may be considered.

\section[short title] {Unlabeled Model}
\setcounter{equation}{0}

Let $${\overline{\nabla}}_{\infty}:=\{x=(x_1,x_2,\dots):x_1\ge
x_2\ge\cdots\ge 0,\ \sum_{i=1}^{\infty}x_i\le 1\}$$ be the closure
of ${\nabla}_{\infty}$ in the product space $[0,1]^{\infty}$. For
$0\le\alpha<1$ and $\theta>-\alpha$, we extend the two parameter
Poisson-Dirichlet distribution $PD(\alpha,\theta)$ from
${\nabla}_{\infty}$ to ${\overline{\nabla}}_{\infty}$. To simplify
notation, we still use $PD(\alpha,\theta)$ to denote this extended
distribution. Let $a(x)$ be the infinite matrix whose $(i,j)$-th
entry is $x_i(\delta_{ij}-x_j)$. Denote by ${\cal P}$ the algebra
generated by 1, $\varphi_2,\ \varphi_3,\dots,\varphi_m,\dots$, where
$\varphi_m(x)=\sum_{i=1}^{\infty}x_i^m$. We consider the bilinear
form ${\cal A}$ of the form
$${\cal
A}(u,v)=\frac{1}{2}\int_{{\overline{\nabla}}_{\infty}}\langle\nabla
u,a(x)\nabla v\rangle dPD(\alpha,\theta),\ \ u,v\in{\cal P}. $$

\begin{thm}\label{ss} The symmetric bilinear form $({\cal A},{\cal P})$ is closable on
$L^2({\overline{\nabla}}_{\infty};PD(\alpha,\theta))$ and its
closure $({\cal A}, D({\cal A}))$ is a regular Dirichlet form.
\end{thm}

\noindent {\bf Proof}\hspace{0.5cm} Define
\begin{equation}\label{s1}
A=\frac{1}{2}\left\{\sum_{i=1}^{\infty}x_i\frac{\partial^2}{\partial
x_i^2}-\sum_{i,j=1}^{\infty}x_ix_j\frac{ \partial^2}{
\partial x_i
\partial x_j}-\sum_{i=1}^{\infty}(\theta
x_i+\alpha)\frac{\partial}{\partial x_i}\right\}.
\end{equation}
The case of $\alpha=0$ corresponds to $A^0$ defined in (\ref{A0}).
One finds that for any $u,v\in{\cal P}$,
$$A^0(uv)=A^0u\cdot v+A^0v\cdot u+\langle\nabla u,a(x)\nabla
v\rangle.$$Hence
\begin{equation}\label{000}
A(uv)=Au\cdot v+Av\cdot u+\langle\nabla u,a(x)\nabla
v\rangle.\end{equation}

We claim that
\begin{equation}\label{2}
\int_{{\overline{\nabla}}_{\infty}}AudPD(\alpha,\theta)=0,\ \
\forall u\in{\cal P}.
\end{equation}
In fact, let $m_1,\dots,m_k\in\{2,3,\dots\}$ and $k\ge 1$. Then we
obtain by (\ref{s1}), (\ref{A0}) and \cite[(2.13)]{EK81} that
\begin{eqnarray}\label{3333}
A(\varphi_{m_1}\cdots\varphi_{m_k})&=&\sum_{i=1}^k\left[{{m_i}\choose{2}}-\frac{m_i\alpha}{2}\right]\varphi_{m_i-1}\prod_{j\not=
i}\varphi_{m_j}+\sum_{i<j}m_im_j\varphi_{m_i+m_j-1}\prod_{l\not=i,j}\varphi_{m_l}\nonumber\\
&
&-\left\{\sum_{i=1}^k\left[{{m_i}\choose{2}}+\frac{m_i\theta}{2}\right]
+\sum_{i<j}m_im_j\right\}\prod_{i=1}^k\varphi_{m_i}\\
&=&\sum_{i=1}^k\left[{{m_i}\choose{2}}-\frac{m_i\alpha}{2}\right]\varphi_{m_i-1}\prod_{j\not=
i}\varphi_{m_j}+\sum_{i<j}m_im_j\varphi_{m_i+m_j-1}\prod_{l\not=i,j}\varphi_{m_l}\nonumber\\
& &-\frac{1}{2}m(m-1+\theta)\prod_{i=1}^k\varphi_{m_i}.\nonumber
\end{eqnarray}
Denote by $\{n_1,n_2,\dots,n_l\}$ an arbitrary partition of
$\{m_1,m_2,\dots,m_k\}$. That is, each
$n_i=m_{i_1}+\cdots+m_{i_{j_i}}$ for some distinct indexes
$i_1,\dots,i_{j_i}$, and
$\{1,2,\dots,k\}=\cup_{i=1}^l\{i_1,\dots,i_{j_i}\}$. By
Ewens-Pitman's sampling formula, we get
\begin{eqnarray*}
&
&\int_{{\overline{\nabla}}_{\infty}}A(\varphi_{m_1}\dots\varphi_{m_k})dPD(\alpha,\theta)\\
&=&\sum_{n_1,n_2,\dots,n_l}
\left\{\sum_{i=1}^l\frac{n_i}{2}(n_i-1-\alpha)\frac{(-\frac{\theta}{\alpha})(-\frac{\theta}{\alpha}-1)
\cdots(-\frac{\theta}{\alpha}-l+1)}{\theta(\theta+1)\cdots(\theta+m-2)}\right.\\
& &\left.\ \ \ \ \cdot\left(\prod_{j\not=i}(-\alpha)\cdots(-\alpha+n_j-1)\right)(-\alpha)\cdots(-\alpha+n_i-2)\right.\\
&
&\left.-\frac{1}{2}m(m-1+\theta)\frac{(-\frac{\theta}{\alpha})(-\frac{\theta}{\alpha}-1)\cdots(-\frac{\theta}{\alpha}-l+1)}{\theta(\theta+1)\cdots(\theta+m-1)}\right.\\
& &\left.\ \ \ \
\cdot\prod_{j}(-\alpha)\cdots(-\alpha+n_j-1)\right\}\\
&=&0,
\end{eqnarray*}
where the value of the right hand side is obtained by continuity
when $\alpha=0$ or $\theta=0$. Similarly, by Ewens-Pitman's
sampling formula, we can further check that
\begin{equation}\label{3}
\int_{{\overline{\nabla}}_{\infty}}(Au)vdPD(\alpha,\theta)=\int_{{\overline{\nabla}}_{\infty}}(Av)udPD(\alpha,\theta),\
\ \forall u,v\in{\cal P}.
\end{equation}

By (\ref{000}), (\ref{2}) and (\ref{3}), we get
$$
{\cal
A}(u,v)=-\int_{{\overline{\nabla}}_{\infty}}(Au)vdPD(\alpha,\theta),\
\ \forall u,v\in{\cal P}.
$$
Hence the symmetric bilinear form $({\cal A},{\cal P})$ is
closable on $L^2({\overline{\nabla}}_{\infty};PD(\alpha,\theta))$
by (\cite[Proposition 3.3]{MR}). The closure $({\cal A}, D({\cal
A}))$ of $({\cal A},{\cal P})$ is a symmetric closed form. To
prove that $({\cal A}, D({\cal A}))$ is a regular Dirichlet form,
it is enough to show that $({\cal A}, D({\cal A}))$ is a Markovian
form. To this end, we will show that $({\cal A}, D({\cal A}))$ is
the same as the closure $({\cal A},{\overline{\cal B}})$ of
$({\cal A},{\cal B})$ with
$$
{\cal B}:=\{u\in
L^2({\overline{\nabla}}_{\infty};PD(\alpha,\theta)):u=f\circ\pi_k\
{\rm for\ some}\ k, f\in C^{\infty}_0(\mathbf{R^k})\},
$$
where $\pi_k: {\overline{\nabla}}_{\infty}\rightarrow
{\mathbf{R^k}}, (x_1,\dots,x_k,\dots)\rightarrow(x_1,\dots,x_k)$.
Note that $({\cal A},{\cal B})$ is clearly Markovian (cf.
\cite[Page 4]{FOT}) and this property is preserved by its closure
(cf. \cite[Theorem 3.1.1]{FOT}).

Let $m\ge 2$. Then one can show that $\varphi_m\in{\overline{\cal
B}}$ by considering the approximation sequence
$\{\varphi^N_m(x):=\sum_{i=1}^{N}x_i^m\}_{N\in{\mathbf{N}}}$. Thus
${\cal P}\subset{\overline{\cal B}}$. To show that ${\cal
B}\subset D({\cal A})$, we need to show that any
finite-dimensional smooth function of the coordinates
$x_1,x_2,\dots$ belongs to $D({\cal A})$. This can be done by
polynomial approximation and noting the fact that
$$
x_1=\lim_{m\rightarrow\infty}(\varphi_m)^{1/m},\ \
x_2=\lim_{m\rightarrow\infty}(\varphi_m-x_1^m)^{1/m},\dots,
$$
where the convergence takes place pointwise on
${\overline{\nabla}}_{\infty}$.

\hfill\fbox

\vskip 0.5cm

It is worth noting that $PD(\alpha,\theta)$ is the unique
probability measure on ${\overline{\nabla}}_{\infty}$ such that
(\ref{2}) is satisfied. In fact, suppose that $\mu\in{\cal
M}_1({\overline{\nabla}}_{\infty})$ satisfying
$$
\int_{{\overline{\nabla}}_{\infty}}Aud\mu=0,\ \ \forall u\in{\cal
P}.
$$
Note that for any $m \geq 2$
$$
A\varphi_m=A^0\varphi_m-\frac{m\alpha}{2}\varphi_{m-1}
=\left[{{m}\choose{2}}-\frac{m\alpha}{2}\right]\varphi_{m-1}-\left[{{m}\choose{2}}+\frac{m\theta}{2}\right]\varphi_{m}.
$$
The fact of $\int_{{\overline{\nabla}}_{\infty}}A\varphi_md\mu=0$
implies that \begin{eqnarray*}
\int_{{\overline{\nabla}}_{\infty}}\varphi_md\mu&=&\frac{\Gamma(m-\alpha)\Gamma(\theta+1)}{\Gamma(1-\alpha)\Gamma(\theta+m)}\\
&=&\frac{\Gamma(\theta+1)}{\Gamma(\theta+\alpha)\Gamma(1-\alpha)}\int_0^1u^m\frac{(1-u)^{\theta+\alpha+1}}{u^{\alpha+1}}du.
\end{eqnarray*}
Then, we obtain by \cite[(6)]{PY} that
\begin{equation}\label{ind}
\int_{{\overline{\nabla}}_{\infty}}\varphi_md\mu=\int_{{\overline{\nabla}}_{\infty}}\varphi_mdPD(\alpha,\theta),\
\ \forall m\in{\mathbf{N}}.
\end{equation}
Furthermore, we obtain by (\ref{3333}), (\ref{ind}) and induction
that
$$
\int_{{\overline{\nabla}}_{\infty}}ud\mu=\int_{{\overline{\nabla}}_{\infty}}udPD(\alpha,\theta),\
\ \forall u\in{\cal P}.
$$
Since ${\cal P}$ is measure-determining, $\mu=PD(\alpha,\theta)$.

By the theory of Dirichlet forms, there exists an essentially
unique Hunt process $(X,(P_x)_{x\in{\overline{\nabla}}_{\infty}})$
on ${\overline{\nabla}}_{\infty}$ with the stationary distribution
$PD(\alpha,\theta)$ such that $X$ is associated with the Dirichlet
form $({\cal A}, D({\cal A}))$ (cf. \cite[Chapter 7]{FOT}). Note
that $A1=0$. By \cite[Proposition 2.3]{Sch3}, one finds that $X$
is a conservative diffusion process. Denote
$P_{PD(\alpha,\theta)}(\cdot)=\int_{{\overline{\nabla}}_{\infty}}P_x(\cdot){PD(\alpha,\theta)}(dx)$.
Then we have the following proposition.
\begin{pro}\label{conser-ergo}
The process $X$ with initial distribution $PD(\alpha,\theta)$ never
leaves $\nabla_{\infty}$, i.e.,
\begin{equation}\label{xcv}
P_{PD(\alpha,\theta)}(X_t\in\nabla_{\infty},\forall t\ge 0)=1.
\end{equation}
In addition, the process $X$ is ergodic, i.e.,
\begin{equation}\label{EK}
\lim_{t\rightarrow\infty}\left\|T_tf-\int_{{\overline{\nabla}}_{\infty}}fdPD(\alpha,\theta)\right\|_{L^2({\overline{\nabla}}_{\infty};PD(\alpha,\theta))}=0,\
\ \forall f\in L^2({\overline{\nabla}}_{\infty};PD(\alpha,\theta)),
\end{equation}
where $(T_t)_{t\ge 0}$ denotes the semigroup associated with $({\cal
A}, D({\cal A}))$ on
$L^2({\overline{\nabla}}_{\infty};PD(\alpha,\theta))$.
\end{pro}
\noindent {\bf Proof}\hspace{0.5cm} We first prove (\ref{xcv}) by
approximation. For $N\in\mathbf{N}$, denote
$\varphi^N_1(x):=\sum_{i=1}^{N}x_i$,
$x\in{\overline{\nabla}}_{\infty}$. Then
$\lim_{N\rightarrow\infty}\|\varphi^N_1-\varphi_1\|_{L^2({\overline{\nabla}}_{\infty};PD(\alpha,\theta))}=0$.
For $N>M$, we have that
\begin{eqnarray*}
{\cal
A}(\varphi^N_1-\varphi^M_1,\varphi^N_1-\varphi^M_1)&\le&\frac{1}{2}\sum_{i=M+1}^N\int_{{\overline{\nabla}}_{\infty}}x_iPD(\alpha,\theta)(dx)\\
&\rightarrow& 0\ \ {\rm as}\ N,M\rightarrow\infty.
\end{eqnarray*}
Thus $\{\varphi^N_1\}_{N\in{\mathbf{N}}}$ is an ${\cal A}$-Cauchy
sequence such that $\varphi^N_1$ converges to $\varphi_1$ in
$L^2({\overline{\nabla}}_{\infty};PD(\alpha,\theta))$ as
$N\rightarrow\infty$. By \cite[Lemma 5.1.2]{FOT}, one finds that
for any $T>0$,
$$
P_{PD(\alpha,\theta)}\left(\sum_{i=1}^NX_i(t)\ {\rm converges\
uniformly\ on}\ [0,T]\ {\rm as}\ N\rightarrow\infty\right)=1.
$$
Then $P_{PD(\alpha,\theta)}(t\rightarrow\sum_{i=1}^{\infty}X_i(t)\
{\rm is\ continuous})=1$. Since for any fixed $t$,
$P_{PD(\alpha,\theta)}(\sum_{i=1}^{\infty}X_i(t)=1)=PD(\alpha,\theta)\{\sum_{i=1}^{\infty}x_i=1\}=1$,
(\ref{xcv}) holds.

Next we turn to the proof of the ergodicity. In fact, it is enough
to verify (\ref{EK}) by considering the following family of
functions
$${\cal
T}:=\{\varphi_{m_1}\cdots\varphi_{m_k}:m_1,\dots,m_k\in\{2,3,\dots\},k\ge
1\}.$$ Let $f\in{\cal T}$. By (\ref{3333}), there exists a
constant $\lambda>0$ and $g\in{\cal T}$ with ${\rm degree}(g)<{\rm
degree}(f)$ such that $Af=-\lambda f+g$. Then
\begin{equation}\label{EK1}
T_tf=e^{-\lambda t}f+e^{-\lambda t}\int_0^te^{\lambda s}T_sgds,\ \
\forall t\ge 0.
\end{equation}
Taking integration on both sides of (\ref{EK1}), we obtain by the
symmetry of $(T_t)_{t\ge 0}$ that
\begin{equation}\label{EK2}
\int_{{\overline{\nabla}}_{\infty}}fdPD(\alpha,\theta)=e^{-\lambda
t}\int_{{\overline{\nabla}}_{\infty}}fdPD(\alpha,\theta)+e^{-\lambda
t}\int_0^te^{\lambda
s}\left(\int_{{\overline{\nabla}}_{\infty}}gdPD(\alpha,\theta)\right)ds.
\end{equation}
Subtracting (\ref{EK2}) from (\ref{EK1}), we get
\begin{eqnarray*}
&
&\left\|T_tf-\int_{{\overline{\nabla}}_{\infty}}fdPD(\alpha,\theta)\right\|_{L^2({\overline{\nabla}}_{\infty};PD(\alpha,\theta))}\\
&\le&e^{-\lambda
t}\left\|f-\int_{{\overline{\nabla}}_{\infty}}fdPD(\alpha,\theta)\right\|_{L^2({\overline{\nabla}}_{\infty};PD(\alpha,\theta))}\\
& &+e^{-\lambda t}\int_0^te^{\lambda
s}\left\|T_sg-\int_{{\overline{\nabla}}_{\infty}}gdPD(\alpha,\theta)\right\|_{L^2({\overline{\nabla}}_{\infty};PD(\alpha,\theta))}ds.
\end{eqnarray*}
Then we can establish (\ref{EK}) by using induction on the degree
of $f$.

\hfill \fbox

\vskip 0.5cm
\begin{rem} The unlabeled two parameter
infinitely-many-neutral-alleles diffusion model considered in this
section is directly motivated by \cite{P}. In \cite{P}, Petrov
used up/down Markov chains and an approximation method to
construct the model. In this section, we use the theory of
Dirichlet forms to give a completely different construction. Our
construction might be more direct and simpler. More importantly,
our observation that the model is given by the classical gradient
Dirichlet form enables us to use this powerful analytic tool to
generalize various basic properties of the
infinitely-many-neutral-alleles diffusion model from the one
parameter setting to the two parameter setting. The Dirichlet form
constructed here differs from the Dirichlet form associated with
the GEM process in \cite{fenwang07} even on symmetric functions.
\end{rem}

There are many problems about the unlabeled two parameter
infinitely-many-neutral-alleles diffusion model which deserve
further investigation. As applications of Theorem \ref{ss}, we
present below several properties of the model via Dirichlet forms,
including a sample path property, a result on large deviations, and
the construction of models with selection.

\begin{thm}\label{st}
Let $X$ be the unlabeled two parameter
infinitely-many-neutral-alleles diffusion model and let $k\ge 1$.
Denote $A_k:=\nabla_{\infty}\cap\{\sum_{i=1}^kx_i=1\}$ and
$D_k:=\nabla_{\infty}\cap\{\sum_{i=1}^kx_i=1\}\cap\{x_k>0\}$.

\noindent (i) If $\theta+\alpha k<1$, then any subset of the
$(k-1)$-dimensional simplex $A_k$ with non-zero
$(k-1)$-dimensional Lebesgue measure is hit by $X$ with positive
probability.

\noindent (ii) If $\theta+\alpha k\ge 1$, then $D_k$ is not hit by
$X$.
\end{thm}

\noindent {\bf Proof}\hspace{0.5cm} We will establish (i) and (ii)
by generalizing \cite[Propositions 2 and 3]{Sch2} to the two
parameter setting. The results are based on Fukushima's classical
result (cf. \cite[Theorem 4.2.1]{FOT}), which says that a Borel
set $B$ is hit by $X$ if and only if $B$ has non-zero capacity. We
use ${\rm Cap}(B)$ to denote the capacity of a Borel set $B$ (cf.
\cite[Chapter 2]{FOT}). Recall that $$ {\rm Cap}(B)=\inf_{B\subset
A\atop{A\ {\rm is\ open}}}{\rm Cap}(A)
$$
and $${\rm Cap}(A)=\inf\left\{{\cal
A}(u,u)+\int_{\overline{\nabla}_{\infty}} u^2dPD(\alpha,\theta):
u\in D({\cal A}), u\ge 1\ {\rm on}\ A,\
PD(\alpha,\theta)-a.e.\right\}$$ if $A$ is an open set.

For $k=1$, let $\nu_1$ denote the Dirac measure at
$(1,0,0,\dots)$. For $k\ge 2$, let
$S_{k-1}:=\{x\in{\mathbf{R^{k-1}}}:x_1\ge\cdots\ge x_{k-1}\ge
0,\sum_{i=1}^{k-1}x_i\le 1\}$ be equipped with $(k-1)$-dimensional
Lebesgue measure and let $\nu_k$ denote the measure induced by the
map $\xi:S_{k-1}\rightarrow A_k$,
$\xi(x_1,\dots,x_{k-1})=(x_1,\dots,x_{k-1},1-\sum_{i=1}^{k-1}x_i,0,0,\dots)$.
In order to show that $A_k$ has non-zero capacity if
$\theta+\alpha k<1$, it is enough to show that there is a
dimension-independent constant $c>0$ such that
\begin{equation}\label{90}
\left(\int_{\overline{\nabla}_{\infty}}ud\nu_k\right)^2\le
c\left({\cal
A}(u,u)+\int_{\overline{\nabla}_{\infty}}u^2dPD(\alpha,\theta)\right),\
\ \forall u\in D({\cal A})\cap C({\overline{\nabla}}_{\infty}).
\end{equation}
For $n\ge k$, denote $B_k=S_n\cap\{\sum_{i=1}^kx_i=1\}$ and use
$\mu_n$, $\nu_{kn}$ to denote respectively the image measures of
$PD(\alpha,\theta)$, $\nu_k$ under the projection of
$\nabla_{\infty}$ onto the first $n$ coordinates. Then (\ref{90})
is equivalent to
\begin{equation}\label{909}
\left(\int_{B_k}fd\nu_{kn}\right)^2\le
c\int_{S_n}\left(\frac{1}{2}\langle\nabla f,a\nabla
f\rangle+f^2\right)d\mu_n,\ \ \forall f\in
C^{\infty}_0({\mathbf{R^n}}).
\end{equation}

To prove (\ref{909}), we will make use of a new coordinate system.
Denote $$r=1-\sum_{i=1}^k(x_i-x_{k+1})$$ and $$
S'_{n-k}=\{x\in{\mathbf{R^{n-k}}}:x_1\ge\cdots\ge x_{n-k}\ge 0,
(k+1)x_1+x_2+\cdots+x_{n-k}\le 1\}.
$$ Consider the map
$\phi:S_n\cap(0<r<1)\rightarrow S_{k-1}\times(0,1)\times
S'_{n-k}$,
\begin{eqnarray*}
\phi(x_1,\dots,x_n)&=&(u_1,\dots,u_{k-1},u_k,u_{k+1},\dots,u_n)\\
&=&
\left(\frac{x_1-x_{k+1}}{r},\dots,\frac{x_{k-1}-x_{k+1}}{r},r,\frac{x_{k+1}}{r},\dots,\frac{x_n}{r}\right).
\end{eqnarray*}
$\phi$ is a one-to-one onto map with the inverse
\begin{eqnarray}\label{inverse}
x_1&=&(1-u_k)u_1+u_ku_{k+1},\nonumber\\
&\vdots&\nonumber\\
x_{k-1}&=&(1-u_{k})u_{k-1}+u_ku_{k+1},\nonumber\\
x_{k}&=&(1-u_{k})(1-(u_1+\cdots+u_{k-1}))+u_ku_{k+1},\nonumber\\
x_{k+1}&=&u_{k+1}u_k,\nonumber\\
&\vdots&\nonumber\\
x_n&=&u_nu_k.
\end{eqnarray}
One can check that the Jocobian of $\phi^{-1}$ is
$(1-u_k)^{k-1}u_k^{n-k}$.

Denote by $h$ the density function of $\mu_n$ with respect to
$n$-dimensional Lebesgue measure. By \cite[Theorem 5.5]{H}, we
have that
\begin{equation}\label{k2}
h(x_1,\dots,x_n)=c_{n,\alpha,\theta}\prod_{j=1}^{n}x_j^{-(\alpha+1)}\left(1-\sum_{i=1}^nx_i\right)^{\theta+\alpha
n-1}\rho_{\alpha,\theta+\alpha
n}\left(\frac{1-\sum_{i=1}^nx_i}{x_n}\right),
\end{equation}
where
\begin{eqnarray*}
c_{n,\alpha,\theta}=\prod_{i=1}^n\frac{\Gamma(\theta+1+(i-1)\alpha)}{\Gamma(1-\alpha)\Gamma(\theta+i\alpha)}=\left\{
\begin{array}{ll}
\theta^n, & \alpha=0,\\
\frac{\Gamma(\theta+1)\Gamma(\theta/\alpha+n)\alpha^{n-1}}{\Gamma(\theta+\alpha
n)\Gamma(\theta/\alpha+1)\Gamma(1-\alpha)^n}, & 0<\alpha<1
\end{array}
\right.
\end{eqnarray*}
and $\rho_{\alpha,\theta+\alpha n}$ is a two parameter version of
Dickman's function, i.e.,
$$
\rho_{\alpha,\theta+\alpha n}(s)=P(sV^{\alpha,\theta+\alpha
n}_1<1),\ \ s\ge 0.
$$
Note that $(1-\sum_{i=1}^nx_i)/x_n$ is only a function of
$u_{k+1},\dots,u_n$ by (\ref{inverse}). Hence we obtain by
(\ref{k2}) that the joint density of $(u_1,\dots,u_n)$ under
$\mu_n\circ\phi^{-1}$ is given by $$
h(u_1,\dots,u_n)=\psi(u_{k+1},\dots,u_n)(x_1\cdots
x_k)^{-(\alpha+1)}(1-u_k)^{k-1}u_k^{\theta+\alpha k-1}
$$
for some function $\psi$. Note that the product $x_1\cdots x_k$ is
only a function of $u_1,\dots,u_{k+1}$. Hence the conditional
density satisfies
$$
h(u_1,\dots,u_k|u_{k+1},\dots,u_n)=\frac{(1-u_k)^{k-1}u_k^{\theta+\alpha
k-1}(x_1\cdots x_k)^{-(\alpha+1)}}{\int_0^1\cdots\int_0^1
(1-u_k)^{k-1}u_k^{\theta+\alpha k-1}(x_1\cdots
x_k)^{-(\alpha+1)}du_1\cdots du_k}.
$$
Therefore, there exists a constant $c_1>0$, which depends on
$\alpha$, $\theta$, $k$ and $\varepsilon>0$ but not $n$, such that
on $\{u_{k+1}\ge\varepsilon\}$ we have
$$
h(u_1,\dots,u_k|u_{k+1},\dots,u_n)\ge
c_1(1-u_k)^{k-1}u_k^{\theta+\alpha k-1}.
$$

We now prove (\ref{909}). Without loss of generality, we assume
that $f$ vanishes for $r(=u_k)\ge 1/2$. In fact, if this condition
is not satisfied, we may obtain (\ref{909}) by multiplying $f$ by
a finite-dimensional smooth function $\gamma\in{\cal B}$, which is
equal to 1 when $r=0$ and vanishes for $r\ge 1/2$. Denote by
$\sigma(du_{k+1},\dots,du_n)$ the distribution of
$u_{k+1},\dots,u_n$ under $\mu_n\circ\phi^{-1}$. We choose
$\varepsilon>0$ such that $p:=\sigma(u_{k+1}>\varepsilon)>0$.
Define $A=S_{k-1}\times(0,1)\times S'_{n-k}$ and
$A_{\varepsilon}=S_{k-1}\times(0,1)\times[S'_{n-k}\cap(u_{k+1}>\varepsilon)]$.
To simplify notation, we denote by $I$ the integral on the left
hand side of (\ref{909}). Then
\begin{eqnarray*}
|I|&=&\left|\int_{S_{k-1}}f\left(u_1,u_2,\dots,u_{k-1},1-\sum_{i=1}^{k-1}u_i,0,\dots,0\right)du_1\cdots
du_{k-1}\right|\\
&=&\left|-\int_0^1\int_{S_{k-1}}\partial_kf((1-u_k)u_1+u_ku_{k+1},\dots,u_nu_k)du_1\cdots
du_{k-1}du_k\right|\\
&=&\left|\frac{1}{p}\int_{A_{\varepsilon}}\partial_k(f\circ\phi^{-1}(u))du_1\dots
du_k\sigma(du_{k+1}\cdots du_n)\right|\\
&\le&\frac{1}{p}\left(\int_{A_{\varepsilon}}u_k[\partial_k(f\circ\phi^{-1}(u))]^2u_k^{\theta+\alpha
k-1}du_1\dots du_k\sigma(du_{k+1}\cdots du_n)\right)^{1/2}\\
& &\times \left(\int_{A_{\varepsilon}}u_k^{-(\theta+\alpha k)}du_1\dots du_k\sigma(du_{k+1}\cdots du_n)\right)^{1/2}\\
&\le&C\left(\int_Au_k[\partial_k(f\circ\phi^{-1}(u))]^2h(u_1,\dots,u_k|u_{k+1},\dots,u_n)du_1\dots du_k\sigma(du_{k+1}\cdots du_n)\right)^{1/2}\\
&\le&C\left(\int_A\langle\nabla f(\phi^{-1}(u)),a(\phi^{-1}(u))\nabla f(\phi^{-1}(u))\rangle\mu_n\circ\phi^{-1}(du)\right)^{1/2}\\
&=&C\int_{S_n}\langle\nabla f,a(x)\nabla f\rangle\mu_n(dx)\\
&=&C{\cal A}(f,f),
\end{eqnarray*}
which proves (\ref{909}). Here $C$ denotes a generic constant whose
value may change from line to line but independent of $n$. For the
last inequality we have used the following estimate
$$
u_k[\partial_k(f\circ\phi^{-1}(u))]^2\le C\langle\nabla
f(\phi^{-1}(u)),a(\phi^{-1}(u))\nabla f(\phi^{-1}(u))\rangle\ {\rm
for}\ u_k\le \frac{1}{2},
$$
which is given by \cite[Lemma 3]{Sch2}.

We now establish (ii). For $k=1$, by (\ref{k2}), we have that
$$
h(x_1)\le cx_1^{-(\alpha+1)}(1-x_1)^{\theta+\alpha-1}
$$
for some constant $c>0$. For $n\ge 1$, we choose $g_n\in
C^{\infty}({\mathbf{R}})$ satisfying $g_n(x)=0$ if $x\le n$,
$g_n(x)=1$ if $x\ge 2n$, and $0\le g_n(x)\le 1$ for all
$x\in{\mathbf{R}}$. Also, we require that $g'_n(x)\le 2/n$ for all
$x\in{\mathbf{R}}$. Set $u_n=g_n\circ\ln((1-x_1)^{-1})$. Then, if
$\theta+\alpha\ge 1$, we have that
\begin{eqnarray*}
{\rm Cap}(A_1)&\le&{\cal
A}(u_n,u_n)+\int_{\overline{\nabla}_{\infty}}
u_n^2dPD(\alpha,\theta)\\
&\le&\frac{c}{2}\int_{1-\exp{(-n)}}^{1-\exp{(-2n)}}[g'_n(\ln(1-x_1)^{-1}]^2(1-x_1)^{-2}x_1(1-x_1)\\
& &\cdot
x_1^{-(\alpha+1)}(1-x_1)^{\theta+\alpha-1}dx_1+PD(\alpha,\theta)\{x_1\ge(1-\exp{(-n)})\}\\
&\le&\frac{2^{1+\alpha}c}{n^2}\int_{1/2}^{1-\exp{(-2n)}}(1-x_1)^{\theta+\alpha-2}dx_1+PD(\alpha,\theta)\{x_1\ge(1-\exp{(-n)})\}\\
&\rightarrow&0\ \ {\rm as}\ n\rightarrow\infty.
\end{eqnarray*}

For $k\ge 2$, we fix an $\varepsilon>0$. Choose $w\in
C^{\infty}({\mathbf{R}})$ satisfying $w(x)=0$ if $x\le\varepsilon$
and $w(x)=1$ if $x>2\varepsilon$. Let $s=\sum_{i=1}^{k}x_i$ and
define $u_n=g_n\circ\ln((1-s)^{-1})$. Set $v_n(x)=u_n(x)w(x_{k})$.
Note that $v_n=1$ on an open subset containing $(s=1)\cap(x_{k}\ge
2\varepsilon)$ and $v_n$ vanishes if $x_{k}\le\varepsilon$. For a
large $n$, the support of $v_n$ is contained in the set
$(1-s)x_k^{-1}\le 1$. Moreover, we obtain by (\ref{k2}) that there
exists a constant $C(\varepsilon, \alpha,\theta)>0$ such that
\begin{equation}\label{mnmn}
h(x_1,\dots,x_k)\le C(\varepsilon,
\alpha,\theta)(1-s)^{\theta+\alpha k-1}\ {\rm on\ the\ support\
of}\ v_n.
\end{equation}
Since $$
\nabla(u_nw)=w\nabla u_n+u_n\nabla w,
$$
we get
\begin{equation}\label{llkk}
{\cal
A}(v_n,v_n)\le\int_{\overline{\nabla}_{\infty}}w^2\langle\nabla
u_n,a\nabla u_n\rangle
dPD(\alpha,\theta)+\int_{\overline{\nabla}_{\infty}}u^2_n\langle
w,a\nabla w\rangle dPD(\alpha,\theta).
\end{equation}
Similar to the $k=1$ case, we can use (\ref{mnmn}) to show that
the first term of the right hand side of (\ref{llkk}) tends to 0
as $n\rightarrow\infty$ if $\theta+\alpha k\ge 1$. Since
$u^2_n\langle w,a\nabla w\rangle\rightarrow 0$ as
$n\rightarrow\infty$, $PD(\alpha,\theta)$-a.e., we conclude that
${\rm Cap}((s=1)\cap(x_k\ge 2\varepsilon))=0$. Since
$\varepsilon>0$ is arbitrary, ${\rm Cap}(D_k)=0$. The proof is
complete.\hfill\fbox

\begin{rem}
In \cite{Sch2}, Schmuland showed that, in the one parameter model,
$A_k$ is hit by $X$ if and only if $\theta<1$. The phase transition
is between infinite and any finite alleles and occurs at $\theta=1$.
In the two parameter model, our Theorem \ref{st} shows that the
phase transition is between infinite and certain finite alleles
(number of alleles is no more than $k_c=[\frac{1-\theta}{\alpha}]$).
The maximum number of finite alleles can be hit is
$[\frac{1}{\alpha}]$ corresponding to $\theta =0$. So the number of
alleles is either infinity or less than or equal to
$[\frac{1}{\alpha}]$. This creates a barrier between finite alleles
and infinite alleles. The results indicate an essential difference
between the one parameter model and the two parameter model, which
deserves a better explanation in terms of coalescent.
\end{rem}

We denote by $(L,D(L))$ the generator of the Dirichlet form
$({\cal A}, D({\cal A}))$ (cf. Theorem \ref{ss}) on
$L^2({\overline{\nabla}}_{\infty}; PD(\alpha,\theta))$. Note that
$Lu=Au$ for all $u\in{\cal P}$, where $A$ is defined in
(\ref{s1}). For $m\ge 2$, define $\lambda_m=m(m-1+\theta)/2$ and
denote by $\pi(m)$ the number of partitions of the integer $m$.
\begin{pro}\label{234} The spectrum of $(L,D(L))$ consists of the eigenvalues
$\{0,-\lambda_2,-\lambda_3,\dots\}$. 0 is a simple eigenvalue and
for each $m\ge 2$, the multiplicity of $-\lambda_m$ is
$\pi(m)-\pi(m-1)$.
\end{pro}

\noindent {\bf Proof}\hspace{0.5cm} The spectrum characterization
has been obtained in \cite{P} using the up/down Markov chains and
approximation. However, a bit more transparent derivation can be
given using our (\ref{3333}). Note that (\ref{3333}) is a
consequence of Pitman's sampling formula and already indicates the
structure of the spectrum of $(L,D(L))$. With \cite[(1.4)]{E}
replaced with our (\ref{3333}), Proposition \ref{234} then follows
from an argument similar to that used in the proof of
\cite[Theorem 2.3]{E}.

\hfill\fbox

\vskip 0.5cm

We now present a result on the large deviations for occupation
time process. It shows that the Dirichlet form $({\cal A}, D({\cal
A}))$ appears naturally as the function governing the large
deviations. Define
$$
L_t(C):=\frac{1}{t}\int_0^t1_C(X_s)ds,\ \ \forall C\in{\cal
B}({\overline{\nabla}}_{\infty}),
$$
where ${\cal B}({\overline{\nabla}}_{\infty})$ denotes the Borel
$\sigma$-algebra of ${\overline{\nabla}}_{\infty}$. We equip
${\cal M}_1({\overline{\nabla}}_{\infty})$ with the
$\tau$-topology, which is generated by open sets of the form
$$
U(\nu;\varepsilon,F):=\left\{\left.\mu\in{\cal
M}_1({\overline{\nabla}}_{\infty})\right|\left|\int Fd\mu-\int
Fd\nu\right|<\varepsilon\right\},
$$
where $\varepsilon>0$, $\nu\in{\cal
M}({\overline{\nabla}}_{\infty})$ and $F\in{
B}_b({\overline{\nabla}}_{\infty})$, the set of bounded Borel
measurable functions on ${\overline{\nabla}}_{\infty}$. The next
result follows from \cite[Theorems 1 and 2]{Mu}.

\begin{pro}\label{236} Let $U$ be a $\tau$-open subset and $K$ be a
$\tau$-compact subset of ${\cal
M}_1({\overline{\nabla}}_{\infty})$. Then for ${\cal A}$-q.e.
$x\in{\overline{\nabla}}_{\infty}$ we have that
$$
\liminf_{t\rightarrow\infty}\frac{1}{t}\log P_x[L_t\in
U]\ge-\inf\{{\cal A}(u,u)|u\in D({\cal A}),u^2PD(\alpha,\theta)\in
U\} $$ and \begin{eqnarray*} & &\inf\left.\left\{\sup_{x\in
{\overline{\nabla}}_{\infty}\backslash
N}\limsup_{t\rightarrow\infty}\frac{1}{t}\log P_x[L_t\in
K]\right|N\subset {\overline{\nabla}}_{\infty}, N\ {\rm is}\ {\cal
A}-{\rm
exceptional}\right\}\\
& &\ \ \ \ \ \ \ \ \ \le-\inf\{{\cal A}(u,u)|u\in D({\cal
A}),u^2PD(\alpha,\theta)\in K\}.
\end{eqnarray*}
\end{pro}

Finally, we would like to point out that the
infinitely-many-neutral-alleles diffusion model considered in this
section can be easily extended to include interactive selection.
\begin{pro}\label{237} Let $\rho\in
L^2({\overline{\nabla}}_{\infty};PD(\alpha,\theta))$ satisfying
$\rho^2\ge\varepsilon>0$, $PD(\alpha,\theta)$-a.e., or $\varphi\in
D({\cal A})$ and $\rho>0$, $PD(\alpha,\theta)$-a.e. Then the
perturbed bilinear form
$$
{\cal
A}^{\varphi}(u,v)=\frac{1}{2}\int_{{\overline{\nabla}}_{\infty}}\langle\nabla
u,a(x)\nabla v\rangle\rho^2dPD(\alpha,\theta),\ \ u,v\in{\cal P}
$$
is closable on $
L^2({\overline{\nabla}}_{\infty};\rho^2PD(\alpha,\theta))$ and its
closure $({\cal A}^{\rho}, D({\cal A}^{\rho}))$ is a regular local
Dirichlet form.
\end{pro}
\noindent {\bf Proof}\hspace{0.5cm} First, we consider the case
that $\rho^2\ge\varepsilon>0$, $PD(\alpha,\theta)$-a.e. Let
$\{u_n\in{\cal P}\}_{n\in{\mathbf{N}}}$ be a sequence satisfying
$u_n\rightarrow 0$ in $
L^2({\overline{\nabla}}_{\infty};\rho^2PD(\alpha,\theta))$ as
$n\rightarrow\infty$ and ${\cal
A}^{\rho}(u_n-u_m,u_n-u_m)\rightarrow 0$ as
$n,m\rightarrow\infty$. Then the strict positivity of $\rho^2$
implies that $\{u_n\}_{n\in{\mathbf{N}}}$ is an ${\cal A}$-Cauchy
sequence and $u_n\rightarrow 0$ in $
L^2({\overline{\nabla}}_{\infty};PD(\alpha,\theta))$. Hence the
closability of $({\cal A},{\cal P})$ implies that
$$\lim_{n\rightarrow\infty}\int_{{\overline{\nabla}}_{\infty}}\langle\nabla
u_n,\nabla u_n\rangle dPD(\alpha,\theta)=0.$$ Thus
$\lim_{n\rightarrow\infty}{\cal A}^{\rho}(u_n,u_n)=0$ by Fatou's
lemma. Therefore $({\cal A}^{\rho},{\cal P})$ is closable.

Now we consider the case that $\rho\in D({\cal A})$ and $\rho>0$,
$PD(\alpha,\theta)$-a.e. Let $(X, P_{PD(\alpha,\theta)})$ be the
Markov process associated with the Dirichlet form $({\cal
A},D({\cal A}))$. Since $\rho\in D({\cal A})$, it has a
quasi-continuous version (cf. \cite[Theorem 2.1.7]{FOT}), which is
denoted by $\tilde{\rho}$. For $n\in{\mathbf{N}}$, we define
$\tau_{n}:=\inf\{t>0:\tilde{\rho}(X_t)\le 1/n\}$ and
$\tau:=\lim_{n\rightarrow\infty}\tau_{n}$. On $\{t<\tau\}$, we
define
$$
M_t^{[\ln\rho]}:=M^{[\ln(\rho\vee(1/n))]}_t,\ \ {\rm if}\
t\le\tau_n,
$$
where $M_t^{[\eta]}$ denotes the martingale part of the Fukushima
decomposition of the additive functional
$\tilde{\eta}(X_t)-\tilde{\eta}(X_0)$ if $\eta\in D({\cal A})$
(cf. \cite[Theorem 5.2.2]{FOT}). We denote by $X^{\rho}$ the
Girsanov transform of $X$ with the multiplicative functional
$L^{[\rho]}_t:=\exp(M_t^{[\ln\rho]}-\frac{1}{2}\langle
M^{[\ln\rho]}\rangle_t)1_{t<\tau}$, where $\langle\cdot\rangle$
denotes the quadratic variation of a martingale. Then $X^{\rho}$
is associated with a Dirichlet form on $
L^2({\overline{\nabla}}_{\infty};\rho^2PD(\alpha,\theta))$ that
extends $({\cal A}^{\rho},{\cal P})$. Therefore $({\cal
A}^{\rho},{\cal P})$ is closable. It is easy to check that its
closure $({\cal A}^{\rho}, D({\cal A}^{\rho}))$ is a regular local
Dirichlet form. The proof is complete.

\hfill\fbox

\section[short title] {Labeled Model}
\setcounter{equation}{0} \noindent In this section, we will
construct measure-valued processes associated with the
two-parameter Dirichlet process through the study of a general
bilinear from. We are successful in two particular cases (cf.
Theorems \ref{S=2} and \ref{S=k} below).

  Let $S$ be a locally compact, separable metric space and $E:={\cal
M}_1(S)$ be the space of probability measures on the Borel
$\sigma$-algebra ${\cal B}(S)$ in $S$. Following (\ref{diri1}), the
two parameter Dirichlet process $\Pi_{\alpha,\theta,\nu_0}$
satisfies
$$
\Pi_{\alpha,\theta,\nu_0}(A)=P\left(\sum_{i=1}^{\infty}\rho_i\delta_{\xi_i}\in
A\right)
$$
for any $A\in{\cal B}(E)$, the Borel $\sigma$-algebra of $E$. We
denote by $E_P$ the expectation with respect to $P$. Set
$$
{\cal F}:={\rm Span}\{\langle f_1,\mu\rangle\cdots \langle
f_k,\mu\rangle:f_1,\dots,f_k\in C_b(S),k\in\mathbf{N}\}.
$$

Consider the following symmetric bilinear form
\begin{equation}\label{finite}
{\cal E}(u,v)=\frac{1}{2}\int_{E}\langle\nabla u(\mu),\nabla
v(\mu)\rangle_{\mu}\Pi_{\alpha,\theta,\nu_0}(d\mu),\ \ u,v\in{\cal
F}.
\end{equation}
Recall that $\nabla u(\mu)$ is the function $$
x\longrightarrow\frac{\partial
u}{\partial\mu(x)}(\mu)=\lim_{\varepsilon\rightarrow
0+}\frac{u((1-\varepsilon)\mu+\varepsilon\delta_x)-u(\mu)}{\varepsilon}
$$
and $\langle f,g\rangle_{\mu}:=\int fgd\mu-(\int fd\mu)(\int
gd\mu)$. Note that $\Gamma(u,v):=\langle\nabla u(\mu),\nabla
v(\mu)\rangle_{\mu}$ is a square field operator. If $({\cal
E},{\cal F})$ is closable on $L^2(E;\Pi_{\alpha,\theta,\nu_0})$,
then following the argument of (\cite[Lemma 7.5 and Proposition
5.11]{Sch}), one can show that the closure $({\cal E},D({\cal
E}))$ of $({\cal E},{\cal F})$ is a quasi-regular local Dirichlet
form. Hence, there exists an essentially unique diffusion process
$X$ which is associated with $({\cal E},D({\cal E}))$ (cf.
\cite[Theorems IV.6.4 and V.1.11]{MR}). This diffusion process is
called the labeled two parameter infinitely-many-neutral-alleles
diffusion model. However, quite different from the unlabeled case,
we find that the closability problem of $({\cal E},{\cal F})$ is
challenging. To understand this point, let us consider the case
that the type space $S$ is finite. This is equivalent to
projecting every $\mu$ in ${\cal M}_1(S)$ to $\{\mu(J_i):
i=1,2,\ldots,n\}$ for certain finite partition $\{J_i:
i=1,2,\ldots,n\}$ of space $S$.

 Let $\{\sigma(t): t \geq 0, \sigma_0=0\}$ be a subordinator with L\'evy measure $x^{-(1+\alpha)}e^{-x}d\,x$, $x>0$,
  and $\{\gamma(t): t \geq 0,
       \gamma_0=0\}$ be a gamma subordinator that is independent of
       $\{\sigma_t: t \geq 0, \sigma_0=0\}$ and has L\'evy measure $x^{-1}e^{-x}d\,x$,
       $x>0$. The next result follows from \cite[Proposition 21]{PY} and the
construction  outlined
      on \cite[Page 254]{Pitman96}.
      \begin{pro}
      {\rm (Pitman and Yor)} Let
      $$
      \gamma(\alpha,\theta)=\frac{\alpha \gamma(\frac{\theta}{\alpha})}{\Gamma(1-\alpha)}.
      $$
      For each $n \geq 1$, and each partition ${J_i: i=1,\ldots,n}$ of $S$, let
      \[
       a_i =\nu_0(J_i),\ \ i=1,\ldots,n,
      \]
      and
     \[
      Z_{\alpha,\theta}(t)=\sigma(\gamma(\alpha,\theta)t),\ \ t \geq 0.
      \]
      Then the distribution of $(\Xi_{\alpha,\theta,\nu_0}(J_1),...,\Xi_{\alpha,\theta,\nu_0}(J_{n}))$
      is the same as the distribution of
      $$\left(\frac{Z_{\alpha,\theta}(a_1)}{Z_{\alpha,\theta}(1)},\dots,\frac{Z_{\alpha,\theta}(\sum_{j=1}^{n}a_j )
      -Z_{\alpha,\theta}(\sum_{j=1}^{n-1}a_j)}{Z_{\alpha,\theta}(1)}\right).$$
      \end{pro}

In general, the distribution function of
$(\Xi_{\alpha,\theta,\nu_0}(J_1),...,\Xi_{\alpha,\theta,\nu_0}(J_{n}))$
cannot be explicitly identified. The exception is the case that
$|S|=2$, i.e., $S$ contains only two elements.

\begin{thm}\label{S=2}
Suppose that $|S|=2$. Then $({\cal E},{\cal F})$ is closable on
$L^2(E;\Pi_{\alpha,\theta,\nu_0})$. Moreover, its closure $({\cal
E},D({\cal E}))$ is a regular local Dirichlet form, which is
associated with a diffusion process on $E$.
\end{thm}

\noindent {\bf Proof}\hspace{0.5cm} We assume without loss of
generality that $0<\alpha<1$ and $\theta>-\alpha$. It is enough to
show that $({\cal E},{\cal F})$ is closable on
$L^2(E;\Pi_{\alpha,\theta,\nu_0})$. Once this is established, the
proof of the last assertion of the theorem is easy. Set
$S=\{1,2\}$, $E=[0,1]$ and $p=1-\bar{p}=\nu_0(1)$. Denote by $dx$
the Lebesgue measure on $[0,1]$. Then ${\cal F}$ is the set of all
polynomials restricted to $[0,1]$ and
$$
{\cal
E}(u,v)=\frac{1}{2}\int_0^1x(1-x)u'(x)v'(x)\Pi_{\alpha,\theta,\nu_0}(dx),\
\ u,v\in{\cal F}.
$$

First, we consider the case that $\theta=0$. It is known (cf.
\cite{lam58}) that
$$
\Pi_{\alpha,0,\nu_0}(dx)=q_{\alpha,0}(x)dx
$$
with
$$q_{\alpha,0}(x)=\frac{p\bar{p}\sin(\alpha\pi)x^{\alpha-1}(1-x)^{\alpha-1}}{\pi[\bar{p}^2x^{2\alpha}+p^2(1-x)^{2\alpha}+2p\bar{p}x^{\alpha}(1-x)^{\alpha}\cos(\alpha\pi)]},\ \ 0\le x\le 1.$$
Define
\begin{eqnarray*}
Lu(x)&=&\frac{1}{2}x(1-x)u^{''}(x)+\frac{\alpha}{2}u'(x)\left[(1-2x)\right.\\
&
&-\left.\frac{2\bar{p}^2x^{2\alpha}(1-x)-2p^2(1-x)^{2\alpha}x+2p\bar{p}(1-2x)x^{\alpha}(1-x)^{\alpha}\cos(\alpha\pi)}
{\bar{p}^2x^{2\alpha}+p^2(1-x)^{2\alpha}+2p\bar{p}x^{\alpha}(1-x)^{\alpha}\cos(\alpha\pi)}\right].
\end{eqnarray*}
Then one can check that $Lu\in L^2(E;\Pi_{\alpha,0,\nu_0})$ for
any $u\in{\cal F}$ and
$$
{\cal E}(u,v)=-\int_0^1(Lu)vd\Pi_{\alpha,0,\nu_0},\ \ u,v\in{\cal
F}.
$$
Therefore $({\cal E},{\cal F})$ is closable on
$L^2(E;\Pi_{\alpha,0,\nu_0})$ by (\cite[Proposition 3.3]{MR}).

We now consider the case that $\theta>0$. To this end, we need to
use a recent result of James et al. By \cite[Example 5.1]{James}
(cf. also \cite[Theorems 3.1 and 5.3]{James}), we have that
$$
\Pi_{\alpha,\theta,\nu_0}(dx)=q_{\alpha,\theta}(x)dx,\ \
q_{\alpha,\theta}(x)=\theta\int_0^x(x-t)^{\theta-1}{\tilde\Delta}_{\alpha,\theta+1}(t)dt.
$$
Here $$
{\tilde\Delta}_{\alpha,\theta+1}(t)=\frac{\gamma_{\alpha-1}(t)\sin(\rho_{\alpha,\theta}(t))-\zeta_{\alpha-1}(t)\cos(\rho_{\alpha,\theta}(t))}{\pi[\zeta_{\alpha}^2(t)
+\gamma_{\alpha}^2(t)]^{(\theta+\alpha)/2\alpha}}
$$
with
$$
\gamma_{d}(t)=\cos(d\pi)t^d\bar{p}+(1-t)^dp,\ \
\zeta_d(t)=\sin(d\pi)t^d\bar{p},\ \ d>-1
$$
and
$$
\rho_{\alpha,\theta}(t)=\frac{\theta}{\alpha}{\rm
arctan}\frac{\zeta_{\alpha}(t)}{\gamma_{\alpha}(t)}+\frac{\pi\theta}{\alpha}1_{\Gamma_{\alpha}}(t),\
\ \Gamma_{\alpha}=\{t\in{\mathbf{R^+}}:\gamma_{\alpha}(t)<0\}.
$$
When $\theta>1$, the expression above can be rewritten as
$$
q_{\alpha,\theta}(x)=(\theta-1)\int_0^x(x-t)^{\theta-2}\Delta_{\alpha,\theta}(t)dt
$$
with
$$
\Delta_{\alpha,\theta}(t)=\frac{\sin\left(\frac{\theta}{\alpha}{\rm
arctan}\left(\frac{\bar{p}\sin(\alpha\pi)t^{\alpha}}{\bar{p}\cos(\alpha\pi)t^{\alpha}+p(1-t)^{\alpha}}\right)+\frac{\pi\theta}{\alpha}
1_{\Gamma_{\alpha}}(t)\right)}
{\pi\{\bar{p}^2t^{2\alpha}+p^2(1-t)^{2\alpha}+2\bar{p}p\cos(\alpha\pi)t^{\alpha}(1-t)^{\alpha}\}^{\theta/2\alpha}},
$$
where $\Gamma_{\alpha}=\emptyset$ if $\alpha\in(0,1/2]$, whereas
$\Gamma_{\alpha}=(0,v_{\alpha}/(1+v_{\alpha}))$ with
$v_{\alpha}=(-p/(\bar{p}\cos(\alpha\pi)))^{1/\alpha}$ if
$\alpha\in(1/2,1)$.

Define
$$
Lu(x)=\frac{1}{2}x(1-x)u^{''}(x)+\frac{1}{2}u'(x)[(1-2x)+x(1-x)q'_{\alpha,\theta}(x)/q_{\alpha,\theta}(x)].
$$
Then one can check that $Lu\in L^2(E;\Pi_{\alpha,\theta,\nu_0})$
for any $u\in{\cal F}$ and
$$
{\cal E}(u,v)=-\int_0^1(Lu)vd\Pi_{\alpha,\theta,\nu_0},\ \
u,v\in{\cal F}.
$$
Therefore $({\cal E},{\cal F})$ is closable on
$L^2(E;\Pi_{\alpha,\theta,\nu_0})$. The proof is
complete.\hfill\fbox

\vskip 0.5cm

From Theorem \ref{S=2}, one can see that even for the one-dimension
case, the generator of the labeled two parameter
infinitely-many-neutral-alleles diffusion model is very complicated.
This indicates an essential difference between the unlabeled model
and the labeled model. More importantly, it explains why it is so
difficult to construct the labeled two parameter
infinitely-many-neutral-alleles diffusion model only using the
ordinary methods that are successful for the one parameter case. So
far we have not been able to solve the closability problem for the
general case. In what follows, we will give further results on the
blinear form $({\cal E},{\cal F})$ and hope they can shed some light
on the problem.

Set
$${\cal
G}:=\{G(\mu)=g(\langle f_1,\mu\rangle,\cdots, \langle
f_k,\mu\rangle),g\in C^{\infty}_b({\mathbf{R^k}}),
f_1,\dots,f_k\in C_b(S)\}.
$$
Let $f\in C_b(S)$ satisfying $\nu_0(f)=0$. We introduce the linear
functional $B_f:{\cal G}\rightarrow \mathbf{R}$ defined by
\begin{equation}\label{22}
B_f(G)=\sum_{s=1}^{\infty}\int
G\left(\sum_{i=1}^{\infty}\rho_i\delta_{\xi_i}\right)f(\xi_s)dP,\
\ G\in{\cal G}.
\end{equation}
Note that (\ref{22}) is well-defined since for $0<\alpha<1$ and
$G(\mu)=g(\langle f_1,\mu\rangle,\cdots, \langle
f_k,\mu\rangle)\in{\cal G}$, we have the following estimate:
\begin{eqnarray*}
\left|\int
G\left(\sum_{i=1}^{\infty}\rho_i\delta_{\xi_i}\right)f(\xi_s)dP\right|&=&\left|\int
\left\{G\left(\sum_{i=1}^{\infty}\rho_i\delta_{\xi_i}\right)-G\left(\sum_{i\not=s}^{\infty}\rho_i\delta_{\xi_i}\right)\right\}f(\xi_s)dP\right|\\
&\le&(\|\partial_1g\cdot f_1\|_{\infty}+\cdots+\|\partial_kg\cdot
f_k\|_{\infty})\|f\|_{\infty}\int\rho_sdP\\
&\le&\frac{c}{s^{1/\alpha}}
\end{eqnarray*}
by \cite[(50)]{PY}, where $c>0$ is a constant which is independent
of $s$.
\begin{pro}\label{lem1} Let $f\in C_b(S)$. Then, for any $v(\mu)=\langle g_1,\mu\rangle\cdots \langle
g_l,\mu\rangle$ with $g_1,\dots,g_l\in C_b(S)$, we have that
$$
{\cal E}(\langle f,\mu\rangle,v)=\frac{\theta}{2}\int_{E}\langle
f-\nu_0(f),\mu\rangle\cdot
v\Pi_{\alpha,\theta,\nu_0}(d\mu)+\frac{\alpha}{2}B_{f-\nu_0(f)}(v).
$$
\end{pro}

\noindent {\bf Proof}\hspace{0.5cm} Let $f\in C_b(S)$ and
$v(\mu)=\langle g_1,\mu\rangle\cdots \langle g_l,\mu\rangle$ with
$g_1,\dots,g_l\in C_b(S)$. Without loss of generality we assume
that $\nu_0(f)=0$. Then
\begin{eqnarray}\label{1}
{\cal E}(\langle f,\mu\rangle,v)&=&\frac{1}{2}\int_{E}\langle
f,\nabla
v(\mu)\rangle_{\mu}\Pi_{\alpha,\theta,\nu_0}(d\mu)\nonumber\\
&=&\frac{1}{2}\int_{E}\sum_{i=1}^l(\langle
fg_i,\mu\rangle\prod_{j\not=i}\langle g_j,\mu\rangle)\Pi_{\alpha,\theta,\nu_0}(d\mu)-\frac{l}{2}\int_{E}\langle f,\mu\rangle\prod_{j=1}^l\langle g_j,\mu\rangle\Pi_{\alpha,\theta,\nu_0}(d\mu)\nonumber\\
&=&\frac{\theta}{2}\int_{E}\langle
f,\mu\rangle\prod_{j=1}^l\langle
g_j,\mu\rangle\Pi_{\alpha,\theta,\nu_0}(d\mu)+\left\{\frac{1}{2}\int_{E}\sum_{i=1}^l(\langle
fg_i,\mu\rangle\prod_{j\not=i}\langle
g_j,\mu\rangle)\Pi_{\alpha,\theta,\nu_0}(d\mu)\right.\nonumber\\
& &-\left.\frac{\theta+l}{2}\int_{E}\langle
f,\mu\rangle\prod_{j=1}^l\langle
g_j,\mu\rangle\Pi_{\alpha,\theta,\nu_0}(d\mu)\right\}.
\end{eqnarray}

Set
$${\cal
H}:=\{\varphi(\mu)=\langle g,\mu^l\rangle:g\in
C_b(S^l),l\in\mathbf{N}\}.
$$
For $l\in\mathbf{N}$, let $\beta=(\beta_1,\beta_2,\dots,\beta_n)$
be an unordered partition of the set $\{1,2,\dots,l\}$. We
associate each $w$ with $\beta_w$ boxes, $1\le w\le n$. Assign the
integers $1,2,\dots,l$ to the $l$ boxes, each box containing
exactly one integer. We denote such an arrangement by $A$. Two
arrangements are said to be the same if they have the same
partition $\beta=(\beta_1,\beta_2,\dots,\beta_n)$ and each $w$,
$1\le w\le n$, is assigned the same (unordered) set of integers.
Define a map $\tau:\{1,2,\dots,l\}\rightarrow\{1,2,\dots,n\}$ by
$\tau(j)=w$ if $j$ is assigned to $w$. Then, we introduce a linear
functional $C_f:{\cal H}\rightarrow \mathbf{R}$ defined by
\begin{eqnarray}\label{10}
&&C_f(\langle g,\mu^l\rangle)\nonumber\\
&&\ \ \ \ =\sum_{{\rm distinct}\ A}\
\left(\frac{(-\frac{\theta}{\alpha})
(-\frac{\theta}{\alpha}-1)\cdots(-\frac{\theta}{\alpha}-(n-1))
\prod_{w=1}^n(-\alpha)(1-\alpha)\cdots(\beta_w-1-\alpha)}{\theta(\theta+1)\cdots(\theta+l-1)}\right.\ \ \ \ \ \ \ \ \ \\
&&\left.\ \ \ \ \ \ \ \ \ \ \ \ \ \ \ \ \ \ \ \cdot\int_{S^n}
g(x_{\tau(1)},\dots,x_{\tau(l)})\sum_{s=1}^nf(x_s)\nu^n_0(dx_1\times\cdots\times
dx_n)\right),\nonumber
\end{eqnarray}
where the value of the right hand side is obtained by continuity
when $\alpha=0$ or $\theta=0$. By (\ref{1}), Pitman's sampling
formula, and comparing the arrangements for sizes $l$ and $l+1$,
we find that
$$
{\cal E}(\langle f,\mu\rangle,v)=\frac{\theta}{2}\int_{E}\langle
f,\mu\rangle\cdot
v\Pi_{\alpha,\theta,\nu_0}(d\mu)+\frac{\alpha}{2}C_f(v).
$$

Let $g\in C_b(S^l)$. Then by (\ref{10}), the assumption that
$\nu_0(f)=0$ and the dominated convergence theorem, we get
\begin{eqnarray}\label{B and C}
C_f(\langle g,\mu^l\rangle)&=&E_P\left\{\sum_{{\rm distinct}\
(i_1,i_2,\dots,i_l)}\left(\rho_{i_1}\rho_{i_2}\dots\rho_{i_l}g(\xi_{i_1},\xi_{i_2},\dots,\xi_{i_l})\sum_{{\rm distinct}\ s\in\{i_1,i_2,\dots,i_l\}}f(\xi_s)\right)\right\}\nonumber\\
&=&\int\sum_{{\rm distinct}\
(i_1,i_2,\dots,i_l)}\left(\rho_{i_1}\rho_{i_2}\dots\rho_{i_l}\sum_{s=1}^{\infty}\int
g(\xi_{i_1},\xi_{i_2},\dots,\xi_{i_l})f(\xi_s)dP^{\xi}\right)dP^{\rho}\nonumber\\
&=&\sum_{s=1}^{\infty}\int\sum_{{\rm distinct}\
(i_1,i_2,\dots,i_l)}\left(\rho_{i_1}\rho_{i_2}\dots\rho_{i_l}\int
g(\xi_{i_1},\xi_{i_2},\dots,\xi_{i_l})f(\xi_s)dP^{\xi}\right)dP^{\rho}\nonumber\\
&=&\sum_{s=1}^{\infty}\int\left\langle
g,\left(\sum_{i=1}^{\infty}\rho_i\delta_{\xi_i}\right)^l\right\rangle
f(\xi_s)dP\nonumber\\
&=&B_f(\langle g,\mu^l\rangle),
\end{eqnarray}
where $P^{\xi}$ and $P^{\rho}$ denote the marginal distributions
of $P$ with respect to $\xi$ and $\rho$, respectively. The proof
is complete.\hfill\fbox

\begin{rem}\label{Re1} If one can show that the linear functional $B_{f-\nu_0(f)}$ defined
by (\ref{22}) is bounded, then there exists a unique $b_f\in
L^2(E;\Pi_{\alpha,\theta,\nu_0})$ such that
$$
B_{f-\nu_0(f)}(G)=\int_Eb_f\cdot Gd\Pi_{\alpha,\theta,\nu_0},\ \
\forall G\in{\cal G}.
$$

We define $$L(\langle f,\cdot\rangle)=-\frac{\theta}{2}\langle
f,\cdot\rangle-\frac{\alpha}{2}b_f.$$ Then
$$
{\cal E}(\langle f,\cdot\rangle,v)=-\int_E(L(\langle
f,\cdot\rangle))vd\Pi_{\alpha,\theta,\nu_0},\ \ \forall v\in{\cal
F}.
$$
In general, we define the operator $L:{\cal F}\rightarrow
L^2(E;\Pi_{\alpha,\theta,\nu_0})$ by induction as follows.
\begin{eqnarray*}
L\left(\prod_{i=1}^k\langle
f_i,\cdot\rangle\right)&=&L\left(\prod_{i=1}^{k-1}\langle
f_i,\cdot\rangle\right)\cdot\langle
f_k,\cdot\rangle+L\left(\langle
f_k,\cdot\rangle\right)\cdot\prod_{i=1}^{k-1}\langle
f_i,\cdot\rangle\\
& &\ \ +\left\langle\nabla\left(\prod_{i=1}^{k-1}\langle
f_i,\cdot\rangle\right),\nabla \langle
f_k,\cdot\rangle\right\rangle.
\end{eqnarray*}
Then one can check that
$$
{\cal E}(u,v)=-\int_E(Lu)vd\Pi_{\alpha,\theta,\nu_0},\ \ \forall
u,v\in{\cal F}.
$$
Therefore, $({\cal E},{\cal F})$ is closable on
$L^2(E;\Pi_{\alpha,\theta,\nu_0})$ by (\cite[Proposition
3.3]{MR}).

If $({\cal E},{\cal F})$ is indeed not closable for the general
case, we may consider its relaxation. We refer the reader to
\cite{Mo} for the definition, existence and uniqueness of
relaxation. The relaxation of $({\cal E},{\cal F})$ is a Dirichlet
form, whose associated Markov process is a good candidate for the
labeled two parameter infinitely-many-neutral-alleles diffusion
model.

\end{rem}

\begin{thm}\label{S=k} Let $S$ be a locally compact, separable metric space, $E={\cal
M}_1(S)$ and $\nu_0\in{\cal M}_1(S)$. Suppose that
$\alpha=-\kappa$ and $\theta=m\kappa$ for some $\kappa>0$ and
$m\in\{2,3,\dots\}$. We denote by $\Pi_{\alpha,\theta,\nu_0}$ the
finite Poisson-Dirichlet distribution. Then the symmetric bilinear
form (\ref{finite}) $({\cal E},{\cal F})$ is closable on
$L^2(E;\Pi_{\alpha,\theta,\nu_0})$. Moreover, its closure $({\cal
E},D({\cal E}))$ is a quasi-regular local Dirichlet form, which is
associated with a diffusion process on $E$.
\end{thm}

\noindent {\bf Proof}\hspace{0.5cm} By independence of the random
variables $\{\xi_s,s=1,2,\dots\}$, we find that
$$\int G\left(\sum_{i=1}^{\infty}\rho_i\delta_{\xi_i}\right)
(f-\nu_0(f))(\xi_s)dP=0,\ \ \forall G\in{\cal G}\ {\rm and}\
s>m.$$ Then the linear functional $B_{f-\nu_0(f)}$ defined by
(\ref{22}) is bounded. Therefore we conclude by Remark \ref{Re1}
that $({\cal E},{\cal F})$ is closable on
$L^2(E;\Pi_{\alpha,\theta,\nu_0})$. Following the argument of
(\cite[Lemma 7.5 and Proposition 5.11]{Sch}), we can further show
that the closure $({\cal E},D({\cal E}))$ of $({\cal E},{\cal F})$
is a quasi-regular local Dirichlet form, which is thus associated
with a diffusion process on $E$. The proof is complete.\hfill\fbox

\vskip 0.5cm

Finally, we present an auxiliary result (cf. Proposition \ref{aux}
below). This result indicates some difficulty of showing the
boundedness of the linear functional $C_f$ defined in (\ref{10}).
Note that the relation between $C_f$ and $B_f$ is described by
(\ref{B and C}). In order to establish the boundedness of $B_f$
and consequently the closability of $({\cal E},{\cal F})$, a
better understanding of the two parameter Poisson-Dirichlet
distributions seems to be needed.

Let $\lambda$ be a partition, i.e., a sequence of the form
$$
\lambda=(\lambda_1,\lambda_2,\dots,\lambda_{l(\lambda)},0,0,\dots),\
\ \lambda_{1}\ge\lambda_2\ge\cdots\ge\lambda_{l(\lambda)}>0,
$$
where $\lambda_i\in \mathbf{N}$.  Denote
$|\lambda|:=\lambda_1+\cdots+\lambda_{l(\lambda)}$. We identify
partitions with Young diagrams. For $k\in\mathbf{N}$, we denote by
$[\lambda:k]$ the number of rows in $\lambda$ of length $k$. For
$n\in\mathbf{N}$, we set (cf. \cite[Page 5]{P})
\begin{eqnarray*}
M_n(\lambda)&:=&\frac{n!}{\prod_{k=1}^{\infty}[\lambda:k]!\cdot\prod_{i=1}^{l(\lambda)}\lambda_i!}\\
& &\cdot\frac{(-\frac{\theta}{\alpha})
(-\frac{\theta}{\alpha}-1)\cdots(-\frac{\theta}{\alpha}-(l(\lambda)-1))\prod_{i=1}^{l(\lambda)}(-\alpha)(1-\alpha)\cdots(\lambda_i-1-\alpha)}{\theta(\theta+1)\cdots(\theta+n-1)}.
\end{eqnarray*}

\begin{pro}\label{aux} Let $0<\alpha<1$ and
$\theta>-\alpha$. Then
\begin{equation}\label{100}
\sum_{\lambda:|\lambda|=n}M_n(\lambda)l(\lambda)=O(n^{\alpha}).
\end{equation}

\noindent {\bf Proof}\hspace{0.5cm} We fix an $n\in\mathbf{N}$.
Let $u(\mu)=\langle 1,\mu\rangle$, $v(\mu)=\langle
1,\mu\rangle\cdots \langle 1,\mu\rangle$ ($n$-fold products),
$f\equiv1$ and $g_1=\cdots=g_n\equiv1$. By considering (\ref{1})
and (\ref{10}), we get
\begin{eqnarray*}
0&=&{\cal E}(u,v)\nonumber\\
&=&\frac{\theta}{2}+\left\{\frac{1}{2}\int_{E}\sum_{i=1}^n(\langle
fg_i,\mu\rangle\prod_{j\not=i}\langle
g_j,\mu\rangle)\Pi_{\alpha,\theta,\nu_0}(d\mu)\right.\nonumber\\
& &-\left.\frac{\theta+n}{2}\int_{E}\langle
f,\mu\rangle\prod_{j=1}^n\langle
g_j,\mu\rangle\Pi_{\alpha,\theta,\nu_0}(d\mu)\right\}\nonumber\\
&=&\frac{\theta}{2}+\frac{\alpha}{2}\sum_{\lambda:|\lambda|=n}M_n(\lambda)l(\lambda)-\frac{\theta+n}{2}\int\sum_{i=1}^{\infty}\rho_i(1-\rho_i)^n\Pi_{\alpha,\theta,\nu_0}(d\mu).
\end{eqnarray*}
Thus, to prove the desired inequality (\ref{100}), we only need to
show that
$$
\sup_{n\in\mathbf{N}}\left\{n^{1-\alpha}\int\sum_{i=1}^{\infty}\rho_i(1-\rho_i)^n\Pi_{\alpha,\theta,\nu_0}(d\mu)\right\}<\infty.
$$

By \cite[(6)]{PY}, we get
\begin{eqnarray*}
n^{1-\alpha}\int\sum_{i=1}^{\infty}\rho_i(1-\rho_i)^n\Pi_{\alpha,\theta,\nu_0}(d\mu)&=&C_1(\alpha,\theta)n^{1-\alpha}\int_0^1u^{-\alpha}(1-u)^{\alpha+\theta+n-1}du\\
&=&C_1(\alpha,\theta)n^{1-\alpha}\cdot{\rm Beta}(1-\alpha,\alpha+\theta+n)\\
&\approx&C_1(\alpha,\theta)n^{1-\alpha}\cdot\Gamma(1-\alpha)(1+\theta+n)^{-(1-\alpha)}\\
&\le&C_2(\alpha,\theta),
\end{eqnarray*}
where $C_1(\alpha,\theta)>0$ and $C_2(\alpha,\theta)>0$ are
constants depending only on $\alpha$ and $\theta$. The proof is
complete.\hfill\fbox

\end{pro}

\vskip 0.5cm
\flushleft {\sc\small S. Feng\\
Department of Mathematics \& Statistics\\
McMaster University\\
1280 Main Street West\\
Hamilton, L8S 4K1, Canada\\
E-mail: shuifeng@mcmaster.ca}

\flushleft {\sc\small W. Sun\\
Department of Mathematics and Statistics\\
Concordia University\\
Montreal, H3G 1M8, Canada\\
E-mail: wsun@mathstat.concordia.ca}

\end{document}